\documentclass[11pt]{article}
\usepackage{a4wide}
\usepackage{latexsym}
\usepackage{amssymb,amsbsy,amsmath,amsfonts,amssymb,amscd}
\usepackage{epsfig, graphicx}
\newcommand{\fer}[1]{(\ref{#1})}

\newcommand{\R}{\mathbb{R}}
\newcommand{\N}{\mathbb{N}}

\newcommand {\al} {\alpha}
\newcommand {\e}  {\varepsilon}

\newcommand {\vp} {\varphi}
\newcommand {\lb} {\lambda}

\newcommand {\Chi} {{\bf \raise 2pt \hbox{$\chi$}} }

\newcommand {\sgn} { {\rm sgn} }

\newcommand {\f}   {\frac}
\newcommand {\p}   {\partial}
\newcommand{\dis}{\displaystyle}
\newcommand {\proof} {\noindent {\bf Proof}. }
\newcommand{\beq}{\begin{equation}}
\newcommand{\beqa}{\begin{eqnarray}}
\newcommand{\bea} {\begin{array}{ll}}
\newcommand{\beqan}{\begin{eqnarray*}}
\newcommand{\eeq}{\end{equation}}
\newcommand{\eeqa}{\end{eqnarray}}
\newcommand{\eeqan}{\end{eqnarray*}}
\newcommand{\eea} {\end{array}}
\newtheorem{theorem}{Theorem}[section]
\newtheorem{lemma}[theorem]{Lemma}

\newtheorem{proposition}[theorem]{Proposition}

\newtheorem{ip}[theorem]{Inverse Problem}
\newcommand{\qed}{{ \hfill
		       {\unskip\kern 6pt\penalty 500
		       \raise -2pt\hbox{\vrule\vbox to 6pt{\hrule width 6pt
		       \vfill\hrule}\vrule} \par}   }}
\title{On the Inverse Problem for a Size-Structured Population Model}

\author{
Beno\^ \i t Perthame\thanks{
D\'epartement de Math\'ematiques et Applications,
\'Ecole Normale Sup\'erieure, CNRS UMR8553 ,
	    45 rue d'Ulm, F~75230 Paris cedex 05, France;
email: perthame@dma.ens.fr}
\thanks{INRIA Rocquencourt, Project BANG, BP 105, F78153 Le Chesnay Cedex}
 \and
Jorge P. Zubelli  \thanks{IMPA, Est. D. Castorina 110, RJ 22460-320, Brazil; email: zubelli@impa.br}
}

\date{\today}

\begin{document}
\maketitle
\pagestyle{plain}
%\tableofcontents
\pagenumbering{arabic}

\begin{abstract}
We consider a size-structured model for cell division and address the question of determining  the division (birth) rate from
the measured stable size distribution of the population. We formulate such question as an inverse problem for
an integro-differential equation posed on the half line. We develop firstly a regular dependency theory for the solution in terms of the coefficients and, secondly, a novel regularization technique for tackling this inverse problem which takes into account the specific nature of the equation. 
Our results rely also on generalized relative entropy estimates and related Poincar\'e inequalities.

\end{abstract}

\noindent {\bf Key-words} Size structure models. Inverse problem. Regularization technique. Entropy methods.
\\
\\
\noindent {\bf AMS Class. Numbers}  35B30, 35R30, 92D25

%%%%%%%%%%%%%%%%%%%%%%%%%%%%%%%%%%%%%%%%%%%%
\section{Introduction}
\label{Sec:intro}
%-------------------------------------------
%%%%%%%%%%%%%%%%%%%%%%%%%%%%%%%%%%%%%%%%%%%%

For many unicellular organisms the mass of the cell, its DNA content, and the level of certain proteins concentrations are often considered as the most relevant parameters for modeling the cell division and, consequently, the population dynamics. Thus, it is natural and usual to consider the evolution of the 
cell density $n$ as a function of the time $t$ and a size
parameter $x$. See \cite{MD,pe, basse} for applications to biology and \cite{bertoin} for other related cases.	In many situations it is of paramount importance to determine from observed data the division (birth) rate $B$ as
a function of the size parameter $x$. However, one is faced with the problem that the cell distribution $n(t,x)$ is not easily measured or
observed. Instead of that, the dynamics of $n(t,x)$ leads to a stable steady distribution $N=N(x)$, which after a suitable re-normalization, can
be observed. In this work we face the problem of determining $B=B(x)$ from noisy and scarce observed data
$N=N(x)$.

Age-structured
populations also constitute basic models for the division of certain cells, 
\cite{chio,asw}. 
There has been substantial interest on the inverse problem for structured population models. 
See for example \cite{EnRunSche,GOP,R89,MD86,PR91} and references therein. These models typically
lead to difficult inverse problems that require a substantial amount of regularization,
especially if one is interested in making use of real data.
In full generality, size structured models have however a somewhat different mathematical structure (only under the 
assumption {\em the smallest daughter cell is larger than half the largest mother cell}, the models lead to similar theory). For this reason,  our approach, differently from that of other works we are aware of, makes use of the
 information given by the so-called {\em stable size distribution} of the model 
 and relies on a number of recent methods based on generalized relative entropy estimates.

As a starting point, we consider the following size-structured	model:
\beq
\label{eq:gencell}
\left \{ \begin{array}{l}
 \f{\p}{\p t} n(t,x) +	\f{\p}{\p x}[g(x) n(t,x)]  + B(x) n(t,x)  = 4  B(2 x) n(t,2 x), \qquad x \geq 0,\, t \geq 0,
\\
\\
n(t,x=0)=0, t> 0,
\\
\\
n(0,x) = n^{0}(x) \ge 0 .
\end{array} \right.
\eeq

This expresses that the evolution of the cell density results from two effects: 
At the one hand, the term $ \f{\p}{\p x} [g(x) n(t,x)] $ which describes the growth of cells 
by nutrient uptake 
with the rate $g(x)$. At the other  hand, the terms containing $B(x)$ and $B(2x)$, which
describe the division of cells of size $2x$ into two cells of size $x$. 
The resulting function-differential equation has a mathematical interest {\em per se}.
One way to understand such interest is that the $x$-derivative leads to transport to
right, whereas the functional term leads to transport to the left. The overall effect
is an equilibrium.
Another way to understand those effects, in mathematical terms, is  to consider the 
evolution of two macroscopic quantities, the total cell number 
$N(t)= \int_0^\infty n(t,x) dx$ and the total biomass 
$M(t)=\int_0^\infty x n(t,x) dx$. Integrating the equation \fer{eq:gencell} yields
\beq \label{eq:cellnum}
\f{d}{dt} N(t)=\int_0^\infty B(x) n(t,x)dx ,
\eeq
this indeed means that number of cells increase only by division, and
\beq \label{eq:cellmass}
\f{d}{dt} M(t)=\int_0^\infty g(x) n(t,x)dx
\eeq
which indeed means that the biomass increases only by nutrient uptake.

For $g\equiv 1$, a case we consider here for simplicity, it was shown in \cite{PR, mmp} that there exists a  unique eigenvalue $\lb_0$ and eigenfunction  $N=N(x)$ such that, after a suitable time re-normalization, the solutions
of (\ref{eq:gencell}) converge to a multiple of $N$ thus  given by the solution of the eigenvalue problem
\beq
\label{eq:celldiv1}
\left \{ \begin{array}{l}
 \f{\p}{\p x} N + (\lb_0+B(x)) N =4  B(2 x) N(2 x), \qquad x \geq 0,
\\
\\
N(x=0)=0,
\\
\\
N(x)>0 \; \text{ for } x>0, \qquad \int_0^\infty N(x)dx =1.
\end{array} \right.
\eeq
Namely, we have under fairly general conditions on the coefficients
$$
n(t,x) e^{-\lb_0 t} {\;}_{\overrightarrow{\; t \rightarrow \infty \; }}\; 
 \rho N(x),
$$
in weighted $L^p$ topologies that are related to entropy properties and that will be described later on. Moreover  exponential rates have been proved to hold for fairly general rates $B$ (\cite{PR,pe}).  Such $N$ is therefore the 
above mentioned {\em stable size distribution}; it is the distribution observed in practice and available for measurements.  
\\

The precise question under consideration in this work is the following:
%  Inverse Problem:
\begin{ip} \label{IP1}
How to recover in a stable way	the birth rate $B(\cdot)$ from noisy data
$N(\cdot)$ and the rate $\lb_0$?
\end{ip}

The plan for this work is the following: In Section~\ref{sec:prelim} we present some
preliminary remarks on a toy model that describes and motivates our approach to the 
regularization of the inverse problem under consideration. 
In Section~\ref{sec:model} we collect a number of results concerning the direct 
problem associated to the model. In particular, concerning the dependence of the
solutions to the model with respect to the coefficients.
In Section~\ref{sec:ipr} we address the inverse problem and its regularization. 
The regularization we propose is novel and consists of introducing the
regularization parameter $\alpha$ directly in the equation as a coefficient
of the perturbed differential equation.
We establish a result concerning the strong stability of such perturbed 
equation and another one demonstrating the consistency of such perturbation
when $\alpha$ goes to zero.
In Section~\ref{sec:estinv} we obtain a convergence rate for this regularized inverse
problem.
Some of the proofs of the results in Section~\ref{sec:model} are rather technical and we 
post-pone them to the Appendix~\ref{sec:ques}.

%%%%%%%%%%%%%%%%%%%%%%%%%%%%%%%%%%%%%%%%%%%%
\section{Preliminaries: A Classical Example Revised} 
\label{sec:prelim}
%-------------------------------------------
%%%%%%%%%%%%%%%%%%%%%%%%%%%%%%%%%%%%%%%%%%%%

In this section we present a simplified problem which shares some 
similarity with \fer{eq:celldiv1}, and allows us to present our strategy for attacking the Inverse Problem \ref{IP1}.  This problem is
the regularization of a classical ill-posed inverse problem, namely recovering
a function from its antiderivative, from a slightly different perspective.
We believe this will make our approach and our estimates more clear
by looking at a specific and simpler well known example.
We refer the reader to the books \cite{BaLe,EnHaNe} for further information and the
classical treatment of the subject.
\\

Inverse problems are characterized by the fact that they are typically ill-posed. Perhaps one of the best well-known
examples is that of differentiation.
Namely, find a function $u$ from its antiderivative $v$.
We focus on the following:
\begin{ip} \label{ip2}
Find  $u: I \rightarrow \mathbb{R}$ such that
\begin{equation} \label{I1}
\int_{0}^{x} \lambda(s) u(s) ds = v (x) \mbox{ , } x \in I
\end{equation}
where $I\subset \mathbb{R}$ is an interval (finite or infinite) of the line with $0$ its left endpoint and
$\lambda : I \rightarrow \mathbb{R}_{>0} $ is given weight function.
\end{ip}
Equation (\ref{I1}) is formally equivalent to
\begin{equation}\label{I2}
\lambda u = \f{\p}{\p x} v \mbox{.}
\end{equation}

We consider the problem~(\ref{I1}), and notice that the map that sends $v \longmapsto u$ is not continuous, when we endow the space of the $v$'s with say the $L^{2}(I ;dx)$ norm and the space of the $u$'s with  $L^{2}(I ;\lambda(x) dx)$ norm, the choice of spaces being motivated by the specific extansion we have in mind. In practice, this reflects itself on the fact that small changes in the measurement of $v$ may lead to huge changes in $u$. Still, one is forced to face such problem in a number of applications.
We now describe some of the main ingredients from inverse problem theory that are used to handle such difficulties as they will be applied to the context of this article. They consist of the following ingredients:
\begin{enumerate}
\item Assume some kind of ``a priori'' bound on the solution $u$ to problem (\ref{I1}), or equivalently on the acceptable data in a sufficiently strong norm. % This is usually referred to as a ``source condition''.
\item Assume that the measured data $v^{\epsilon}$ satisfies $||v - v^{\epsilon}|| \le \epsilon $ in a suitable weak norm $|| \cdot ||$.
\item Approximate the problem (\ref{I1}) by one that is well-posed for the weak norm and solve the stabilized problem instead of the unstable
one choosing the regularization parameter in an optimal way.
\end{enumerate}
Let us illustrate the technique in the problem under consideration. For more information we refer the reader to \cite{BaLe,EnHaNe}.
Let us consider
\begin{equation} \label{Ialpha}
\left\{
\begin{array}{l}
\alpha \f{\p}{\p x}(\lambda u_\alpha) + \lambda u_\alpha  = \f{\p}{\p x} v \mbox{ ,} \\
\lambda u_{\alpha} (0) = 0
\end{array}
\right.
\end{equation}
In what follows we shall assume that $v(0) = 0 $ whenever $v$ is defined in a neighborhood of $0$ due to some extra assumptions (such as $\f{\p}{\p x} v \in L^2 $).

In order to implement the above program to the problem under consideration, we shall establish the following two inequalities:
%------------------------------
\begin{lemma}
If $u_{\alpha}$ satisfies (\ref{Ialpha}), then
\begin{equation} \label{star}
|| u_{\alpha} ||_{L^{2}(I, \lambda^2 dx) } \le \frac{1}{\alpha} || v ||_{L^{2}(I, dx) } \mbox{ .}
\end{equation}
If $\f{\p^2}{\p x^2 } v \in L^{2}(I,dx) $ and $u_\alpha$ satisfies (\ref{Ialpha}), then
\begin{equation} \label{starstar}
|| u_{\alpha} ||_{L^{2}(I, \lambda^2 dx) } \le \alpha ||\f{\p^2}{\p x^2 } v ||_{L^{2}(I, dx) } \mbox{ .}
\end{equation}
\end{lemma}
%---------------------------------------------
As a direct consequence of Equation~(\ref{starstar}) we have that if $u_{\alpha}$ is a solution of equation (\ref{Ialpha}) and
$\overline{u}$ is a solution of (\ref{I2}), then the difference $u_{\alpha} - \overline{u}$ satisfies
\begin{equation} \label{star2}
|| u_{\alpha} - \overline{u} ||_{L^{2}(I, \lambda^2 dx) } \le \alpha ||\f{\p^2}{\p x^2 } v ||_{L^{2}(I, dx) } \mbox{ .}
\end{equation}
We now join the above estimates with the ``a priori'' assumption that
\begin{equation} \label{apriori}
\big|\big| \f{\p^2}{\p x^2 } v \big|\big|_{L^2(I,dx)} \le E \mbox{ ,}
\end{equation}
and consider a problem with noisy data given by $v^{\epsilon} $ such that
\begin{equation} \label{dataeps}
|| v - v^{\epsilon} ||_{L^2(I, dx)} \le \epsilon  \mbox{ .}
\end{equation}
Thus, upon solving the regularized problem (\ref{Ialpha}) with data $v^{\epsilon}$ we are left with a
solution $u_{\alpha}^{\epsilon} $ of
\begin{equation} \label{Ialpha2}
\left\{
\begin{array}{l}
\alpha \f{\p}{\p x} (\lambda u_\alpha^{\epsilon})  + \lambda u_\alpha^{\epsilon}  = \f{\p}{\p x} v^{\epsilon} \mbox{ ,} \\
\lambda u_{\alpha}^{\epsilon} (0) = 0 \mbox{ .}
\end{array}
\right.
\end{equation}
We further assume that $\f{\p}{\p x} v^{\epsilon} \in L^{2}(I,dx)$.
Now we consider the distance between $u_\alpha^{\epsilon}$ and $\overline{u}$
\begin{eqnarray*}
|| u_\alpha^\epsilon - \overline{u} ||_{L^2(I, \lambda^2 dx) } & \le & || u_\alpha - \overline{u} ||_{L^2(I, \lambda^2 dx) } + || u_{\alpha}^{\epsilon} -  u_\alpha ||_{L^{2}(I, \lambda^2 dx) } \\
 & \le & \alpha ||\f{\p^2}{\p x^2 } v ||_{L^{2}(I, dx) } + \frac{1}{\alpha} || v^{\epsilon} - v ||_{L^{2}(I, dx) } \\
 & \le & \alpha E + \frac{\epsilon}{\alpha} \; .
\end{eqnarray*}
Now, we are free to choose $\alpha$ in such a way as to minimize the R.H.S. of the last inequality.
Using the elementary fact that $\alpha \mapsto (\alpha E + \epsilon/\alpha) $ has a minimum at $\alpha = \sqrt{E/\epsilon}$
and choosing such optimal $\alpha = \alpha(\epsilon) $ we get that
\begin{equation} \label{fine}
|| u_\alpha^{\epsilon} - \overline{u} ||_{L^{2}(I, \lambda^2 dx) } \le 2 \epsilon^{1/2} E^{1/2} \mbox{ .}
\end{equation}

Summing up, if we know an ``a priori'' bound on the data $v$ and $v^{\epsilon}$ of the form (\ref{apriori}) and the data is
measured with an accuracy $\epsilon$ then we get that the stabilized computed solution is within
$\mathcal{O}(\epsilon^{1/2} E^{1/2}) $	of the actual solution $\overline{u}$ in the norm of the space
$ L^{2}(I, \lambda^2 dx)$.

The reasoning developed above can be pushed in different directions and in the next subsection we will show how it can be extended further.

We remark that Equation~(\ref{Ialpha}) is just one of the different possibilities of regularizing the Problem (\ref{I2}).
Another possibility would be to solve an equation of the form
\begin{equation} \label{IIalpha}
- \alpha \f{\p^2}{\p x^2 } \frac{u_{\alpha}}{\lambda} + \lambda u_{\alpha} = \f{\p}{\p x} v
\end{equation}
subject to appropriate boundary conditions.
It is an easy exercise to verify that Equation~(\ref{IIalpha}) formally corresponds to Tikhonov regularization of the Inverse Problem~\ref{ip2}.
The form of Equation~(\ref{IIalpha}), although very attractive and general is not well adapted to our hyperbolic problem and we do not consider it here.

%%%%%%%%%%%%%%%%%%%%%%%%%%%%%%%%%%%%%%%%%%%%
%-------------------------------------------
\section{Model Properties and Regular Dependency upon the Coefficients}
\label{sec:model}
%-------------------------------------------
%%%%%%%%%%%%%%%%%%%%%%%%%%%%%%%%%%%%%%%%%%%%

In this section we shall address a preliminary issue, namely the
smoothness of {\em direct}  map
$$
B \mapsto (\lb_0,N) \mbox{ }
$$
in suitable function spaces. 
Then, the Inverse Problem \ref{IP1} makes sense in such spaces. 
This requires some assumptions and preliminary mathematical material that we introduce now.
\\

Firstly,
we need some assumptions. Here we make the choice of simplicity, keeping in mind that more realistic hypotheses are possible that would lead to more technical proofs, which we prefer to avoid. 
The measurable division rates $B$ are supposed to satisfy
% Assumption B_m, B_M
\beq
\label{as:1}
\exists B_m,\; B_M \quad \text{such that } \; 0< B_m \leq B(x) \leq B_M,
\eeq

Secondly, because the adjoint  equation is as important as  the direct eigenvalue problem, we consider both the cell division equation, to find $(\lb_0,N)$ solution to
\beq
\label{eq:celldiv}
\left \{ \begin{array}{l}
 \f{\p}{\p x} N + (\lb_0+B(x)) N =4  B(2 x) N(2 x), \qquad x \geq 0,
\\
\\
N(x=0)=0,
\\
\\
N(x)>0 \; \text{ for } x>0, \qquad \int_0^\infty N(x)dx =1.
\end{array} \right.
\eeq
and the adjoint equation
\beq
\label{eq:celldivd}
\left \{ \begin{array}{l}
 \f{\p}{\p x} \vp - (\lb_0+B(x)) \vp =- 2 B(x ) \vp(\frac x 2), \qquad x \geq 0,
\\
\\
\vp(x)>0, \qquad \int_0^\infty \vp(x) N(x)dx =1.
\end{array} \right.
\eeq
The solution $\vp$ is used in several places because it gives natural bounds for the direct problem. See the proof of Theorem \ref{th:regularity} and the entropy dissipation in  Appendix A.1.
\\

The existence and uniqueness of solutions to these eigenproblems was proved in \cite{PR},  and it is known that $N$ as well as its derivative vanish at $x=0$ and $x=\infty$. We also recall the following results that we will use later on.

%\begin{figure}{\centering
%\includegraphics[width =.45\textwidth, height =.45\textwidth]{sizestuc1.pdf}
%\quad 
%\includegraphics[width =.45\textwidth, height =.45\textwidth]{sizestruc2.pdf} 
%\vspace{-.1cm}
%\caption{Solution to the spectral problem  
%\fer{eq:celldiv} with two different birth rates. 
%Left: $B=1$. Right: $B=Min(2, 16 \ x^2)$.} }
%\end{figure}

\begin{theorem}{\cite{PR, michelth,pe}} Under the assumption~\fer{as:1}, the solution $N$ to \fer{eq:celldiv} satisfies the properties
\beq
\label{eq:f1}
 B_m \leq \lb_0 = \int B(x) N(x) dx \leq B_M,
\eeq
\beq
\label{eq:f2}
\f{1}{ B_M} \leq  \int x N(x) dx =\f{1}{\lb_0 } \leq \f{1}{B_m},
\eeq
\beq
\label{eq:f3}
0 \leq N(x) \leq 2 B_M ,
\eeq
\beq
\label{eq:f4}
e^{ax} N \in L^1\cap L^\infty(\R^+),  \qquad \forall  a <\lb_0+B_m, \qquad \vp (x) \leq C(1+x).
\eeq
\label{th:estimates}
\end{theorem}

\proof We only prove the estimates on $N$ and refer to \cite{michelth} for the sub-linearity result on $\vp$ and for more accurate results on the behavior of $N$ at infinity, in particular, when $B$ has a limit $B_\infty$ at infinity, one can obtain exponential decay with $a=\lb_0+B_\infty$.

The inequalities \fer{eq:f1} and \fer{eq:f2} are obtained by integrating Equation	 \fer{eq:celldiv}  against respectively the weights $dx$ and $x \; dx$.

For \fer{eq:f3}, we use	 \fer{eq:f1}, and from the Equation \fer{eq:celldiv}, we obtain
$$
 \f{\p}{\p x} N(x)  \leq 4 B(2x) N(2x) \Longrightarrow N(x) \leq 4 \int_0^x B(2y) N(2y)dy \leq 2B_M.
$$

For  \fer{eq:f4}, we multiply  \fer{eq:celldiv}	 by $e^{ax}$ and integrate between 
$0$ and $\infty$. In fact, a standard truncation argument is needed for complete 
justification. See \cite{pe} for details. We obtain
\beq
\label{eq:estexp}
\begin{array}{rl}
 \int_0^\infty [\lb_0+B(x)-a] e^{ax} N(x) dx &= 4 \int_0^\infty B(2x) N(2x) e^{ax}dx
\\ \\
& \leq	2\int_0^{\infty} B(x) N(x) e^{ax/2}dx.
\end{array}
\eeq
With the choice $a=\lb_0$, since the right hand side is bounded by
$$
 2 \left(\int_0^{\infty} B(x) N(x) e^{ax}dx \right)^{1/2} \left( \int_0^{\infty} B(x) N(x)  \right)^{1/2} = 2 \left(\int_0^{\infty} B(x) N(x) e^{ax}dx \right)^{1/2} \lb_0^{1/2},
$$
we obtain that
\beq
\label{eq:f5}
\int_0^{\infty} B(x) N(x) e^{\lb_0 x}dx \leq 4 \lb_0 \leq 4B_M.
\eeq
Then, we can iterate on $a$ and choose any $a<B_m+ \lb_0 \leq 2 \lb_0$ while keeping the right hand side of \fer{eq:estexp} bounded. This proves the $L^1$ statement of	 \fer{eq:f4}.
\\

It remains to use \fer{eq:celldiv} again to deduce, thanks to the chain rule, that $\f{\p}{\p x} [e^{ax} N(x)] \in L^1(\R^+)$, and the result follows.
\qed

\bigskip

With these results and methods of proof, we can now state the main new result of this section.
%%%%%%%%%%%%%%%%%%%%%%%
\begin{theorem} Under the assumption~\fer{as:1},  the map $B \mapsto (\lb_0,N)$, from $L^\infty(\R^+)$ into $[B_m,B_M] \times L^1\cap L^\infty(\R^+)$ is
\\
(i) \; Continuous in $B$ under the weak-$\ast$ topology of $L^\infty(\R^+)$, 
\\
(ii) \ Locally Lipschitz continuous  in $B$ under the strong topology of $L^2(\R^+)$ into $L^2(\R^+)$, namely \fer{eq:esteig} and \fer{eq:spectr} below hold, for $m$ large enough,
\\
(iii) Class $C^1$ in the spaces of statement (ii). 
\label{th:regularity}
\end{theorem}
%%%%%%%%%%%%%%%%%%%%%%%

In statement (i) the  weak-$\ast$ topology of $L^\infty(\R^+)$ is natural from the assumption \fer{as:1}. For estimates (ii) and (iii), we have chosen the $L^2$ space which is used later on for the inverse problem where it plays an essential role.
\\

\proof First of all, we introduce some notations that we will use throughout the proof. 
For two functions $B$, $\bar B$ satisfying \fer{as:1}, we set
\beq
\label{eq:defdelta}
\delta B =  \bar B- B, \quad \Delta  =\| \bar B- B\|_{L^2 (\R^+)}.
\eeq
 For the corresponding solutions we define the differences
$$
\delta	N= \bar N- N, \quad \delta  \lb= \bar \lb_0- \lb_0.
$$
(i) This first statement follows thanks to  the strong uniform  estimates on $N$  in Theorem \ref{th:estimates}. Indeed, consider  a sequence $B_n$ that
weak-$\ast$ converges to $B$ and the corresponding solutions $N_n$ to \fer{eq:celldiv}.  Theorem \ref{th:estimates} 
provides us with uniform bounds on $N_n$ and thus with uniform bounds on $\f{\p N_n}{\p x}$. This means that the sequence $N_n$ is relatively strongly compact. Therefore, it converges to the unique $N$ solution of Equation \fer{eq:celldiv} for the limit $B$.
\\
\\
(ii) Firstly, we estimate the difference between the eigenvalues. As a consequence of \fer{eq:celldiv} for $\bar N$ (multiplied by $\vp$) and \fer{eq:celldivd} for $\vp$ (multiplied by $\bar N$), we obtain
$$
\delta \lb \int \bar N \; \vp = \int \delta B \vp\big(\f{x}{2}\big) \bar N +  \int \delta B \vp	 \bar N.
$$
Therefore, according to  the {\em a priori} estimates of Theorem \ref{th:estimates} on $\vp$ and on $N$,
\beq
| \delta \lb|  \leq C(B, \bar B)  \Delta,\label{eq:esteig}
\eeq
with $C(B, \bar B)= \|\f{\vp}{1+x}\|_{L^\infty(\R^+)}	\| (1+x) \bar  N \|_{L^2(\R^+)}/\int \bar N \; \vp$.
\\

Secondly,  we estimate $\delta N$. Writing the difference of the solutions to the cell-division equation, we have
\beq
 \f{\p}{\p x} \delta N+	 (\lb_0+B(x)) \;   \delta N= 4 B(2x)  \delta N(2x) + \delta R(x),
\label{eq:eqdn}
\eeq
with
\beq
\delta R(x) =  4 \delta B(2x) \; \bar N(2x) + \big(\delta \lb+\delta B(x)\big) \bar N(x).
\label{eq:eqdr}
\eeq
We have, since $N$ has exponential decay at infinity using again Theorem \ref{th:regularity}, and for all $m\in \N$, 
\beq
\|  \delta R  (1+  x^m) \|_{L^2(\R^+)} \leq  C(B, \bar B) \; \Delta,
\label{eq:estrest}
\eeq
for another constant $C$ still controlled by $B_m$ an d$B_M$. 
Also, since $\vp$ has at most linear growth at infinity, we know that
\beq
\int_0^\infty |	 \delta R(x) |\vp(x) dx  \leq	C(B, \bar B) \; \Delta.
\label{eq:estrestp}
\eeq

It is also useful to recall that the Equation \fer{eq:eqdn} is spectral ($N$ is the Krein-Rutman eigenvector) and one readily checks the solvability and the uniqueness conditions
\beq
\int_{\R^+}  \varphi \; \delta R =0, \qquad \int_{\R^+}	 \delta N=0.
\label{eq:eqdnsolv}
\eeq

Our local lipschitz regularity result then follows from the two next lemmas

%---------------------------------------------------------------
\begin{lemma}[Entropy dissipation]
Assume that Equation~\fer{as:1} holds and that 
$N$ and $\vp$ are the solutions of \fer{eq:celldiv} and \fer{eq:celldivd}. 
Then, the solution $\delta N$ to \fer{eq:eqdn}, with the conditions \fer{eq:eqdnsolv}, 
satisfies for all convex function $H:\R\to \R$,
\beq  \begin{array}{rl}
 \dis{ \int_0^\infty } 4\vp(x) B(2x) N(2x) & \left[ H'\big( \f{\delta N(x)}{N(x)} \big)\big( \f{\delta N(2x)}{N(2x)} - \f{\delta N(x)}{N(x)} \big)
-H \big( \f{\delta N(x)}{N(x)} \big)+ H\big( \f{\delta N(2x)}{N(2x)}\big) \right] dx
\\
\\
&= \dis{ \int_0^\infty} H'\big( \f{\delta N(x)}{N(x)} \big) \; \delta R(x) \; \vp(x) \; dx.
\end{array}
\label{eq:entr1}
\eeq
\label{lm:entropy}
\end{lemma}

This lemma is a variant of the generalized relative entropy inequalities introduced in \cite{MPR,mmp}.  Its proof will be presented in the Appendix~\ref{sec:ques}.
We point out, however, the main difficulty that arises here: Since $N$ and $\bar N$ do not have necessarily the same exponential decay at infinity, one cannot use the natural weighted quantities $\int_0^\infty \vp N |\f{\delta R}{N}|^p$  which are not necessarily finite. Therefore, we are sometimes forced to work in $L^1$. This type of entropy inequality is related to Poincar\'e type inequalities which would be a way to conclude the proof at this stage. This area is still a very active
field of research for so-called hypocoercive operators. See for example \cite{villani}. Here, it does not seem to follow from standard methods and we prefer to prove directly the following consequence of Lemma \fer{lm:entropy}

%-----------------------------------------------------------
\begin{lemma}[Spectral gap] There is a $\nu >0$ such that
\beq
\nu \| \delta N\|_{L^2(\R^+)} \leq \| \delta R \; (1+ x^m) \|_{L^2(\R^+)} \leq C \Delta,
\label{eq:spectr}
\eeq
where the parameter $m$ (large enough) is chosen so that $\lb_0 >   \f{B_M}{2^{m-1}}$.
\label{lm:spectral}
\end{lemma}
%-----------------------------------------------------------
We again refer to the Appendix for a proof.
\\

This Lemma expresses the Lipschitz regularity stated in point (ii).
\\
\\
(iii) The third claim is now a general fact on the Fr\'echet derivatives of mapping with quadratic nonlinearities. Whenever it is Lipschitz continuous on a subset, it is $C^1$ by usual algebraic manipulations. This concludes the proof of Theorem \ref{th:regularity}.
\qed

%%%%%%%%%%%%%%%%%%%%%%%%%%%%%%%%%%%%%%%%%%%%
\section{The Inverse Problem and its Regularization}
\label{sec:ipr}
%-------------------------------------------
%%%%%%%%%%%%%%%%%%%%%%%%%%%%%%%%%%%%%%%%%%%%
%%%%

The inverse problem consists in finding $B$ from the knowledge of the population growth rate and the cell density $(\lb_0,N)$. If one could guarantee that the measurement $N$ is very smooth one could directly consider to solve the	equation on $B$, equivalent to the cell-division Equation \fer{eq:celldiv} written with $y=2x$,
\beq
4B(y) N(y) = B\big(\f y 2\big) N\big(\f y 2)+ \lb_0 N\big(\f y 2\big)+2 \f{\p}{\p y} N\big(\f y 2\big), \qquad y>0.
\label{eq:ex}
\eeq
This is a well-posed equation on $B$ as long as $N$ satisfies regularity properties such as $ \f{\p}{\p y} N\big(\f y 2\big) \in L^p$ for some $p\geq 1$. 
But our interest lies on the fact that fluctuations on $N$ make it a mere $L^p$ function because, given a set of measurements on $N$, one does not have a way of controlling the precision of the measurements on $\f{\p}{\p y}N$. This motivates working with a regularized problem
\beq
\left\{
\begin{array}{l}
\al \f{\p}{\p y} (B_\al N) + 4B_\al(y) N(y) = B_\al\big(\f y 2\big) N\big(\f y 2\big)+ \lb_0 N\big(\f y 2\big)+2 \f{\p}{\p y} N\big(\f y 2\big), \qquad y>0,
\\
\\
B_\al N(0)=0,
\end{array} \right.
\label{eq:invreg}
\eeq
where $0<\al<1$ is a small parameter adapted to the level of noise as explained in Section \ref{sec:prelim} for the toy model. Notice that this is still a well-posed problem in appropriate function
spaces when written with a general source $F(y)$
\beq
\left\{
\begin{array}{l}
\al \f{\p}{\p y} (B_\al N) + 4B_\al(y) N(y) = B_\al\big(\f y 2\big) N\big(\f y 2\big)+ F(y), \qquad y>0,
\\
\\
B_\al N(0)=0,
\end{array} \right.
\label{eq:invregbis}
\eeq

Before we analyze the issue of estimating the inverse problem, let us explain why this is indeed  a well-posed problem and for simplicity we restrict ourselves to $L^2$. In particular, we did not impose boundary conditions because we expect that the point $x=0$ is characteristic since the cell density satisfies
\beq
N\in H^2(\R^+), \qquad N(x) >0 \quad \text{for }\; x>0, \qquad N(0)=0 \qquad \f{\p}{\p x} N(0)=0.
\label{eq:invas1}
\eeq
%-------------------------------
\begin{theorem}[Strong stability] Assume \fer{eq:invas1} and that $F\in L^2$.
 Then, there is a unique solution to \fer{eq:invregbis} such that $BN \in H^1$
with  
$$
\al |B_\al N(y)|^2+  \int |B_\al N|^2 \leq C   \int |F|^2 , \qquad \forall y \geq 0,
$$
$$
\al^2	 \int |\f{\p}{\p y}B_\al N|^2 \leq C   \int |F|^2 ,
$$
and 
$$
 \int |\f{\p}{\p y} B_\al N|^2 \leq C	\int |\f{\p}{\p y}F|^2 \qquad  \text{ if } \; F(0)=0.
$$
\label{th:regest}
\end{theorem}
%-------------------------------

\proof We drop the index $\al$ in order to simplify the notation. This model is a variant of the cell division equation and the Cauchy-Lipschitz theory applies and gives the existence for small $y$. The important point is to establish the  {\em a priori} bounds. For the first set of estimates, we multiply the Equation \fer{eq:invregbis} by $BN$ and integrate from $0$ to $y$. This yields
$$
\begin{array}{rl}
 \f{\al}{2} |BN(y)|^2 + &4 \int_0^y |BN|^2(s)ds
\\
\\
 &=\dis{ \int_0^y (BN)(z) \; (BN)\big(\f z 2\big) dz + \int_0^y (BN)(z) \; F(z) dz}
\\
\\
&\dis{ \leq \f 1 2  \int_0^y |BN|^2+ \f 1 2 \int_0^y |BN\big(\f z 2\big)|^2 dz+\f 1 2  \int_0^y |BN|^2+\f 1 2 \int_0^y	\; |F(z)|^2 dz}.
\end{array}
$$
Altogether, this gives the two bounds claimed  in the first statement above.
\\

The second estimate, on $\al  \int |\f{\p}{\p y}BN|^2$, is then a consequence of the first two using Equation~\fer{eq:invregbis} once again.
\\

For the third inequality, we differentiate the equation and write $Q= \f{\p}{\p y}BN$. We find
$$
\al \f{\p}{\p y}Q + 4 Q(y)= \f 1 2 Q\big(\f y 2\big) + \f{\p}{\p y}F.
$$
Therefore, multiplying by $Q$ and integrating by parts gives
$$
4 \int_0^\infty Q^2 \leq \f{\al}{2} Q^2(0) + \f 1 2 \int_0^\infty Q(y) \;  Q\big(\f y 2\big) dy +\int_0^\infty  Q(y) \f{\p}{\p y}F(y) dy,
$$
and thus
$$
4 \int_0^\infty Q^2 \leq \f{\al}{2} Q^2(0) + \f 1 4 \int_0^\infty Q^2+ \f 1 4 \int_0^\infty Q^2\big(\f y 2\big) dy +
\f 1 2 \int_0^\infty Q^2 +  \f 1 2 \int_0^\infty |\f{\p}{\p y}F|^2.
$$
Altogether this gives
$$
\frac{11}{4}\int_0^\infty Q^2 \leq \f{\al}{2} Q^2(0) +\int_0^\infty |\f{\p}{\p y}F|^2.
$$
It remains to notice that from the Equation \fer{eq:invregbis} on $BN$ and the condition $N(0)=0$ we have also
$\al Q(0)= \al \f{\p}{\p y} BN(0) =F(0)$ and the result is proved.
\qed

\medskip

The first question we want to answer is how much  this regularization on $B$ differs from its exact value
%-------------------------------
\begin{theorem}[Consistency] Assume \fer{eq:invas1}, then the solutions to \fer{eq:ex} and \fer{eq:invreg} satisfy
$$
\int |B_\al-B|^2 N \leq C \al^2 \int \left[ |\f{\p}{\p y} N|^2+ |\f{\p^2}{\p y^2} N|^2\right].
$$
\label{th:consist}
\end{theorem}
%-------------------------------

\proof We write $\delta B = B_\al-B$ and we have
$$
\al \f{\p}{\p y} (\delta B  \; N) + 4\delta B(y) \; N(y) = \delta B \big(\f y 2\big) \; N\big(\f y 2\big)- \al \f{\p}{\p y} ( B	 N)  .
$$
This is again an equation of the form \fer{eq:invregbis} with source term
$$
F= - \al \f{\p}{\p y} ( B  N)
$$
and thus we can apply the first estimate of Theorem \ref{th:regest} and obtain
$$
\int |\delta B	\; N|^2 \leq C \al^2 \int \left| \f{\p}{\p y} ( B  N) \right|^2.
$$
Next we apply the second estimate of Theorem \ref{th:regest}  to \fer{eq:ex} with $F= \lb_0 N+ \f{\p}{\p y} N$ and obtain
$$
 \int \left| \f{\p}{\p y} ( B  N) \right|^2 \leq C \int [ |\f{\p}{\p y} N|^2+ |\f{\p^2}{\p y^2} N|^2].
$$
Altogether these two estimates give the announced inequality.
\qed

%%%%%%%%%%%%%%%%%%%%%%%%%%%%%%%%%%%%%%%%%%%%
\section{Recovering Estimates for the Regularized Inverse Problem}
\label{sec:estinv}
%-------------------------------------------
%%%%%%%%%%%%%%%%%%%%%%%%%%%%%%%%%%%%%%%%%%%%

Consider now a smooth solution $(\lb_0, N)$ to Equation \fer{eq:celldivd} corresponding to division rate $B$, that we want to recover from a noisy measurement	$N_\e$ with
\beq
\int |N- N_\e|^2 \leq \e^2, \qquad N_- \leq N_\e \leq N_+,
\label{eq:regas1}
\eeq
where $N_\pm$ are two smooth functions that can serve as a filter for unrealistic data. We assume their behavior contains possible properties of the true solution and serves to assert the boundary condition is well defined with $N_+(0)=N_+(\infty)=0$. Notice $N_-=0$ is a possibility.
The regularized method is aimed at furnishing an approximation $B_{\e,\al}$ of the exact coefficient $B$ through Equation \fer{eq:invreg} with $N_\e$ in place of $N$,
\beq
\left\{
\begin{array}{l}
\al \f{\p}{\p y} (B_{\e,\al} N_\e) + 4B_{\e,\al}(y) N_\e(y) = B_{\e,\al}\big(\f y 2\big) N_\e\big(\f y 2\big)+ \lb_0 N\big(\f y 2\big)+2 \f{\p}{\p y} N_\e\big(\f y 2\big), \qquad y>0,
\\
\\
(B_{\e,\al} N_\e)(0)=0.
\end{array} \right.
\label{eq:invfull}
\eeq
We prove the following error estimate on the coefficient recovered by the above regularization procedure
%-------------------------------
\begin{theorem}[Convergence rate] Assume \fer{eq:invas1} and \fer{eq:regas1}, then the solutions  $B$ to \fer{eq:ex} and $B_{\e,\al}$ to \fer{eq:invfull},  satisfy the error estimate,
$$
\|B_{\e,\al}-B\|_{L^2( N_\e^2dx)}  \leq C \al \| N\|_{H^2}+  \f{C+\|B\|_{L^\infty}}{\al} \|N_{\e}-N\|_{L^2}.
$$
\label{th:cvrate}
\end{theorem}
%-------------------------------

This theorem relies on a first estimate which expresses weak stability and improves Theorem \ref{th:regest} using the special structure in the RHS of Equation \fer{eq:invfull}.  Namely, we have
%-------------------------------
\begin{proposition}[Weak stability] Assume \fer{eq:invas1} and \fer{eq:regas1}, given two solutions $B_\al$ to \fer{eq:invreg} and \fer{eq:invfull} satisfy
$$
\int |B_{\e,\al}N_{\e} -B_\al N|^2  \leq \f{C}{\al^2} \|N_\e-N\|^2_{L^2}.
$$
\label{pr:cvrate}
\end{proposition}
%-------------------------------

\proof We set $Q=B_{\e,\al}N_{\e,\al} -BN$, $R=N_\e-N$	which satisfies
$$
\al \f{\p}{\p y}Q +4 Q(y) =Q\big(\f y 2\big) +\lb_0 R+ 2\f{\p}{\p y}R\big(\f y 2\big) .
$$
Then, we obtain after multiplication by $Q$ and integration
\begin{eqnarray*}
\f{\al}{2} [Q^2(y)-Q^2(0)] + 4\int_0^y |Q(z)|^2 dz & = & \int_0^y Q(z) \; Q(\f z 2)\;	 dz+ \\
  & & \lb_0 \int_0^y QR-2 \int_0^y R\big(\f z 2\big) \f{\p}{\p y}Q+ 2R\big(\f y 2\big)Q(y)-2QR(0). \end{eqnarray*}
Using again the equation on $Q$ (and the Cauchy-Schwarz inequality) we find
$$
\begin{array}{rl}
\f{\al}{2}[Q^2(y)-Q^2(0)] + 2 \int_0^y |Q(z)|^2 dz\leq &C \int_0^y |R|^2+2Q(y)R\big(\f y 2\big)-2QR(0)
\\
\\
&-\f{2}{\al}\int_0^y R\big(\f z 2\big)	[Q\big(\f z 2\big)  +\lb_0 R(z)+ \\
& \f{\p}{\p y}R\big(\f z 2\big) -4Q(z)]dy .
\end{array}
$$
Therefore
$$
\f{\al}{2} |Q(y)-\f{2}{\al} R\big(\f y 2\big)|^2 + \int_0^y |Q(z)|^2 dz\leq
\f{\al}{2} |Q(0)-\f{2}{\al} R(0)|^2 + \f{C}{\al^2} \int_0^y |R|^2.
$$
Notice that even though $R$ is not defined point-wise, $Q(y)-\f{2}{\al} R\big(\f y 2\big)$ is well defined (as a $H^1$ function) and it vanishes due to \fer{eq:regas1}.
\qed

\bigskip

We are now ready to prove our main result, Theorem \ref{th:cvrate}. We write
$$
\begin{array}{rl}
\|B_{\e,\al}-B\|^2_{L^2( N_\e^2dx)}  &\leq 2\int |B_{\e,\al}N_{\e} -B N|^2+ 2\int |B N_{\e} -B N|^2
\\
\\
&\leq 4\int |B_{\e,\al}N_{\e} -B_\al N|^2+ 4\int |B_\al N -B N|^2+  2\int |B N_{\e} -B N|^2.
\end{array}
$$
The first term is controlled thanks to Proposition \ref{pr:cvrate} and gives the second error term in the estimate of
Theorem \ref{th:cvrate}. The second  term is controlled because of Theorem \ref{th:consist} and gives the first error term. The third term gives the last contribution once $\|B\|_{L^\infty}$ has been factored out and the Theorem  \ref{th:cvrate} is proved.
\qed

%----------------------------------------------------------------
\section{Conclusion and Suggestions for Further Research}
%----------------------------------------------------------------
%----------------------------------------------------------------
%

We have proposed a non-standard regularization method for the Inverse Problem associated with the cell-division equations. It is based on a first order operator that leads us to solving a new equation whose structure is very close to the equation itself (but for a different unknown). This allowed us to develop a consistency and convergence analysis that can be used in practice because the outcome are standard inequalities taking  into account the noise in the measured data. 
\\

Our theoretical approach has been focused on the equal mitosis, when mother  cells divide into two equal daughter cells. A first extension of the method would be to deal with more general division equations as 
$$
\left \{ \begin{array}{l}
 \f{\p}{\p t} n(t,x) +	\f{\p}{\p x}[g(x) n(t,x)]  + B(x) n(t,x)  =2  \int_0^\infty B(y) \beta (\f x y) n(t,y) \f{dy}{y}, \qquad x \geq 0,\, t \geq 0,
\\
\\
n(t,x=0)=0, t> 0.
\end{array} \right.
$$
The function $\beta$ represents now the repartition of daughter cell sizes and a natural inverse problem would be to recover the division rate $B$ from measurements of the cell density $N$ once $\beta$ is known. In principle, our method can extended but technical estimates have to be reformulated in this new context.
\\

The numerical validation of the procedure proposed here is straightforward because it relies on a regularized equation that has a standard form and numerical methods are available. As far as validation on real data is concerned, let us 
point out that present experimental devices allow us to measure not only the size repartitions in a cellular culture, but also the molecular content of certain representative proteins. Mathematical models describing the division process based on molecular contents are being developped. To our knowledge, identification of coefficients in this context is still largely open.

\section*{Acknowledgments} J.P.Z. acknowledges financial support from CNPq, grants 302161/2003-1 and
474085/2003-1.
This research was made possible through the Brazil-France cooperation program, which
is gratefully acknowledged by both authors.

%%%%%%%%%%%%%%%%%%%%%%%%%%%%%%%%%%%%
%-------------------------------------------
\appendix
\section{Appendix}
\label{sec:ques}
%-------------------------------------------
%%%%%%%%%%%%%%%%%%%%%%%%%%%%%%%%%%%%

This appendix is devoted to the proof of two technical lemmas that were used in the Section~\ref{sec:model} to obtain Lipschitz regularity of the solution to the eigenproblem for the cell division equation.

\subsection{Generalized Relative Entropy. Proof of Lemma \ref{lm:entropy}.}
Using \fer{eq:eqdn} and \fer{eq:celldiv},  we have
\beq
\f{\p}{\p x} \f{\delta N(x)}{N(x)}  =
4B(2x) \f{N(2x)}{N(x)} \big[ \f{\delta N(2x)}{N(2x)} -\f{\delta N(x)}{N(x)} \big] + \f{\delta R(x)}{N(x)},
\label{ap:1}
\eeq
and thus, for $H(\cdot)$ a convex function,
$$
  \f{\p}{\p x}H\big( \f{\delta N(x)}{N(x)}  \big) = 4 B(2x) \f{N(2x)}{N(x)} H'\big( \f{\delta N(x)}{N(x)} \big) \big[ \f{ \delta N(2x)}{N(2x)} -\f{\delta N(x)}{N(x)}\big] + H'\big( \f{\delta N(x)}{N(x)} \big) \; \f{\delta R(x)}{N(x)} .
$$
On the other hand, combining Equations \fer{eq:celldiv} and \fer{eq:celldivd}, we have
$$
\f{\p}{\p x} \big( N(x) \vp(x) \big)= 4 \vp(x) B(2x)  N(2x) -2 N(x) B(x) \vp(\f{x}{2}).
$$
Therefore, 
$$
\begin{array}{rl}
  \f{\p}{\p x} \big[ N(x) \vp(x) H\big( \f{\delta N(x)}{N(x)}  \big) \big]
& = 4 B(2x)N(2x) \vp(x) H'\big( \f{\delta N(x)}{N(x)} \big) \big[\f{\delta N(2x)}{N(2x)} - \f{\delta N(x)}{N(x)}\big]
\\
\\
& +\big[4 \vp(x) B(2x)	N(2x) -2 N(x) B(x) \vp(\f{x}{2}) \big] H\big( \f{\delta N(x)}{N(x)} \big)
\\
\\
&+ H'\big( \f{\delta N(x)}{N(x)} \big) \; \delta R(x) \; \vp(x).
\end{array}
$$
After integration in $x$ we arrive at
$$ \begin{array}{rl}
0 &= \dis{ \int_0^\infty } 4\vp(x) B(2x) N(2x)	H'\big( \f{\delta N(x)}{N(x)} \big) \big[\f{\delta N(2x)}{N(2x)} - \f{\delta N(x)}{N(x)}\big] dx
\\
\\
&+ \dis{ \int_0^\infty}4 \vp(x) B(2x)  N(2x) \big[ H \big( \f{\delta N(x)}{N(x)} \big)- H\big( \f{\delta N(2x)}{N(2x)}\big) \big] dx
\\
\\
&+  \dis{ \int_0^\infty} H'\big( \f{\delta N(x)}{N(x)} \big) \; \delta R(x) \; \vp(x) \; dx.
\end{array}
$$
Once the terms are reorganized, this is exactly the statement of Lemma \ref{lm:entropy}.
\qed

\bigskip
%------------------------------------------------------------------------------------------------------

\subsection{A Poincar\'e-Type Inequality. Proof of Lemma \ref{lm:spectral}.}
We argue by contradiction and assume that for a sequence  $\bar B_k \to B$,  the corresponding $\nu_k$ vanishes. Then, there is a family of $\delta R_k$, $\delta N_k$ that satisfy
\beq
\nu_k \| \delta N_k\|_{L^2(\R^+)} \geq \| \delta R_k \; (1+ x^m) \|_{L^2(\R^+)}.
\label{eq:spectrcontr}
\eeq
After re-normalizing (by multiplication), we can always assume that
$$
\| \delta N_k\|_{L^2(\R^+)} =1, \qquad \| \delta R_k(1+x^m)\|_{L^2(\R^+)} {\;}_{\overrightarrow{\; k \rightarrow \infty \; }}\;	 0 .
$$
Then,  our proof consists in several steps:  First, we prove compactness and then pass to the limit, next we use the generalized relative entropy to identify the limiting solution and finally, prove a contradiction.
\\

\noindent  {\it Compactness}. Notice that the Equation \fer{eq:eqdn} automatically implies that (recall $\vp>0$ is smooth and sublinear at infinity), that
$\delta N_k$ is bounded in $H^1(\R^+)$. This provides us with local compactness.  It remains to obtain a control at infinity;  we also notice that
$$
\f 1 2	\f{\p}{\p x} (\delta N_k)^2 +	(\lb_0+B(x)) \;	 (\delta N_k)^2= 4 B(2x) \delta N_k(2x) \; \delta N_k(x) + \delta R_k(x) \; \delta N_k(x) \mbox{ ,}
$$
which, after integration in $x$ with the weight	 $x^m$ gives for all $a>0$
$$
\begin{array}{rl}
- \f m 2  \int_{\R^+}x^{m-1} (\delta N_k)^2  + \int_{\R^+}x^m (\lb_0+B(x)) \;  (\delta N_k)^2&
\leq 4 \int_{\R^+}x^m B(2x) [a \delta N_k(2x)^2 +\f 1 a	 \delta N_k(x)^2]
\\
\\
&+\sqrt{ \int_{\R^+}x^m (\delta N_k)^2 \int_{\R^+}x^m (\delta R_k)^2}.
\end{array}
$$
Therefore,  we also obtain
$$
\begin{array}{rl}
 \int_{\R^+}x^m (\lb_0+B(x)) \;	 (\delta N_k)^2&
\leq  \f m 2  \int_{\R^+}x^{m-1} (\delta N_k)^2	 + a\; 2^{-m+1} \int_{\R^+}x^m B(x)  \delta N_k(x)^2
\\
\\
&+\f{B_M}{a} \int_{\R^+}x^m  \delta N_k(x)^2
+\sqrt{ \int_{\R^+}x^m (\delta N_k)^2 \int_{\R^+}x^m (\delta R_k)^2}.
\end{array}
$$
With $a=2^{m-1}$, we conclude that
$$
\begin{array}{rl}
\lb_0 \int_{\R^+}x^m  (\delta N_k)^2&
\leq  \f m 2  \int_{\R^+}x^{m-1} (\delta N_k)^2	 + \f{B_M}{a} \int_{\R^+}x^m ( \delta N_k)^2
+\sqrt{ \int_{\R^+}x^m (\delta N_k)^2 \int_{\R^+}x^m (\delta R_k)^2}
\\
\\
&\leq C	 \int_{\R^+} (\delta N_k)^2  + (\f{\lb_0}{2}+\f{B_M}{a}) \int_{\R^+}x^m ( \delta N_k)^2 +\sqrt{ \int_{\R^+}x^m (\delta N_k)^2 \int_{\R^+}x^m (\delta R_k)^2}.
\end{array}
$$
It remains to choose $m$ large enough so that $\lb_0> 2 \f{B_M}{a}$ and we conclude on a uniform control
\beq
\label{eq:unifren}
\int_{\R^+}x^m	(\delta N_k)^2 \leq C\left[1+ \int_{\R^+} x^m  (\delta R_k)^2\right],
\eeq
that concludes the global compactness in $L^2(\R^+)$ of the sequence $\delta N_k$.
\\

\noindent {\it Limit}. Therefore, due to the above facts and using Lemma~\ref{lm:entropy}, we may pass to the strong limit and obtain
$$
\delta N_k \to n\in H^1(\R^+) \qquad  x^m n \in L^2(\R^+) , \, \, \forall m>0,
$$
$$
 \f{\p}{\p x} n+  (\lb_0+B(x)) \;  n = 4 B(2x)	n(2x), \quad \int_{\R^+} \delta n=0.
$$
Notice that, thanks to the moment estimate, we also have for all $m>0$,
$$
(1+x^m) n \in L^1(\R^+), \mbox{ and }  \int_{\R^+} |n| \vp < \infty.
$$
This allows us to give a meaning to the vanishing integral. This follows from the Cauchy-Schwarz inequality
$$
\int_{\R^+} (1+x^m) |n| \leq \left( \int_{\R^+} (1+x^m)^2 (1+|x^2|)  |n|^2(x) dx \;	\int_{\R^+} \f{dx}{1+x^2} \right)^{1/2}.
$$
\\

\noindent  {\it Application of the GRE}. We can now choose the family of convex functions $H(u)=(u-\xi)_+$ in the Generalized Relative Entropy of Lemma \ref{lm:entropy}, with $\xi >0$ a parameter. We find, because we now deal with $\delta R=0$, after reorganizing the terms in \fer{eq:entr1}, that
$$
 \dis{ \int_0^\infty } 4\vp(x) B(2x) N(2x) \big( \f{n(2x)}{N(2x)} - \f{n(x)}{N(x)} \big)
\left(
 \big( \sgn_+(\f{n(x)}{N(x)}-\xi \big)	- \sgn_+\big( \f{n(2x)}{N(2x)} -\xi\big) dx =0.
\right)
$$
In fact, one needs to justify it by a preliminary truncation of the integral at $x=0$  and  $x\approx \infty$, we leave to the reader the corresponding analysis. Again, let us point out that in $L^1$ this is possible because the corresponding $x$-derivative term
$$
 [\vp(x) N(x) \big( \f{n(x)}{N(x)} -\xi \big)_+ ]
$$
is  integrable (but for the square entropy we arrive at $[\vp(x) N(x) \big( \f{n(x)}{N(x)} \big)^2 ]$ which  is not integrable!).
\\

\noindent  {\it Conclusion}. From the previous step, we conclude that, for all $x >0$, $\xi>0$, we have
$$
\sgn_+(\f{n(x)}{N(x)}-\xi \big)	 = \sgn_+\big( \f{n(2x)}{N(2x)} -\xi\big) ,
$$
which implies $ \f{n(x)}{N(x)} = \f{n(2x)}{N(2x)} $ and thus, tells us that the limit satisfies
$$
\f{\p}{\p x} \f{n(x)}{N(x)}=0.
$$
From 
 $ \f{n(x)}{N(x)} = \f{n(2x)}{N(2x)} $, we conclude that $n = CN$ and the 
vanishing total integral
of $n$ implies that $n=0$.
This gives the sought contradiction and concludes the proof of Lemma \ref{lm:spectral}.
\qed

%%%%%%%%%%%%%%%%%%%%%%%%%%%%%%%%%%%%%%%%
%%%%%%%%%%%%%%%%%%%%%%%%%%%%%%%%%%%%
%
%%%%%% BIBLIO %%%%%%%%%%%%%%%%%%%%%%%%%%
%
%%%%%%%%%%%%%%%%%%%%%%%%%%%%%%%%%%%%

\bibliographystyle{unsrt}
\bibliography{zupe}

\end{document}